\documentclass{gtmon_a}
\pdfoutput=1
\usepackage{graphicx}

\proceedingstitle{Groups, homotopy and configuration spaces (Tokyo
  2005)}
\conferencestart{5 July 2005}
\conferenceend{11 July 2005}
\conferencename{Groups, homotopy and configuration spaces,
                in honour of Fred Cohen's 60th birthday}
\conferencelocation{University of Tokyo, Japan}

\editor{Norio Iwase}
\givenname{Norio}
\surname{Iwase}

\editor{Toshitake Kohno}
\givenname{Toshitake}
\surname{Kohno}

\editor{Ran Levi}
\givenname{Ran}
\surname{Levi}

\editor{Dai Tamaki}
\givenname{Dai}
\surname{Tamaki}

\editor{Jie Wu}
\givenname{Jie}
\surname{Wu}

\title{Poisson structures on the homology of the space of knots}

\author{Keiichi Sakai}
\givenname{Keiichi}
\surname{Sakai}
\address{Graduate School of Mathematical Science\\University of Tokyo\\\newline
3-8-1 Komaba\\Meguro\\Tokyo 153-8914\\Japan}
\email{ksakai@ms.u-tokyo.ac.jp}
\urladdr{}



\arxivreference{math.AT/0608326}

\volumenumber{13}
\issuenumber{}
\publicationyear{2008}
\papernumber{21}
\startpage{463}
\endpage{482}

\doi{}
\MR{}
\Zbl{}

\received{19 September 2006}
\revised{16 July 2007}
\accepted{23 July 2007}
\published{19 March 2008}
\publishedonline{19 March 2008}
\proposed{}
\seconded{}
\corresponding{}
\version{}

\makeatletter
\def\cnewtheorem#1[#2]#3{\newtheorem{#1}{#3}[section]
\expandafter\let\csname c@#1\endcsname\c@thm}

\let\xysavmatrix\xymatrix
\def\xymatrix{\disablesubscriptcorrection\xysavmatrix}
\AtBeginDocument{\let\bar\wbar\let\tilde\wtilde\let\hat\what}
\newcommand{\we}{\smash{\rlap{\kern 6pt
\raise 4pt\hbox{\footnotesize $\sim$}}}\longrightarrow}
\newcommand{\wel}{\smash{\rlap{\kern 8pt
\raise 4pt\hbox{\footnotesize $\sim$}}}\longleftarrow}

\newtheorem{thm}{Theorem}[section]
\cnewtheorem{lem}[thm]{Lemma}
\cnewtheorem{prop}[thm]{Proposition}
\cnewtheorem{cor}[thm]{Corollary}

\theoremstyle{definition}
\cnewtheorem{rem}[thm]{Remark}

\makeatother  %  move after \newtheorem block

\newcommand{\A}{\mathcal{A}}
\newcommand{\calD}{\mathcal{D}}
\newcommand{\calO}{\mathcal{O}}

\newcommand{\abs}[1]{\lvert {#1} \rvert}
\newcommand{\ang}[2]{\langle {#1},\, {#2}\rangle}

\newcommand{\id}{\mathrm{id}}
\newcommand{\img}{\mathrm{im}\,}
\newcommand{\Imm}{\mathrm{Imm}_n}
\newcommand{\K}{\mathcal{K}_n}
\newcommand{\Map}{\mathrm{Map}\,}
\newcommand{\rank}{\mathrm{rank}\,}
\newcommand{\Tot}{\mathrm{Tot}\,}
\newcommand{\tTot}{\widetilde{\Tot}}
\newcommand{\X}{X^{\bullet}_n}

\begin{document}

\begin{htmlabstract}
In this article we study the Poisson algebra structure on the homology
of the totalization of a fibrant cosimplicial space associated with an
operad with multiplication.  This structure is given as the Browder
operation induced by the action of little disks operad, which was
found by McClure and Smith.  We show that the Browder operation
coincides with the Gerstenhaber bracket on the Hochschild homology,
which appears as the E<sup>2</sup>-term of the homology spectral sequence
constructed by Bousfield.  In particular we consider a variant of the
space of long knots in higher dimensional Euclidean space, and show
that Sinha's homology spectral sequence computes the Poisson algebra
structure of the homology of the space.  The Browder operation
produces a homology class which does not directly correspond to chord
diagrams.
\end{htmlabstract}

\begin{abstract}
In this article we study the Poisson algebra structure on the homology
of the totalization of a fibrant cosimplicial space associated with an
operad with multiplication.  This structure is given as the Browder
operation induced by the action of little disks operad, which was
found by McClure and Smith.  We show that the Browder operation
coincides with the Gerstenhaber bracket on the Hochschild homology,
which appears as the $E^2$-term of the homology spectral sequence
constructed by Bousfield.  In particular we consider a variant of the
space of long knots in higher dimensional Euclidean space, and show
that Sinha's homology spectral sequence computes the Poisson algebra
structure of the homology of the space.  The Browder operation
produces a homology class which does not directly correspond to chord
diagrams.
\end{abstract}

\begin{asciiabstract}
In this article we study the Poisson algebra structure on the homology
of the totalization of a fibrant cosimplicial space associated with an
operad with multiplication. This structure is given as the Browder
operation induced by the action of little disks operad, which was
found by McClure and Smith.  We show that the Browder operation
coincides with the Gerstenhaber bracket on the Hochschild homology,
which appears as the E2-term of the homology spectral sequence
constructed by Bousfield. In particular we consider a variant of the
space of long knots in higher dimensional Euclidean space, and show
that Sinha's homology spectral sequence computes the Poisson algebra
structure of the homology of the space.  The Browder operation
produces a homology class which does not directly correspond to chord
diagrams.
\end{asciiabstract}

\maketitle

%%% 1st section %%%

\section{Introduction}

In \cite{McClureSmith04_2} McClure and Smith proved that the totalization of a cosimplicial space
$\calO^{\bullet}$ associated with a non-symmetric topological operad $\calO$ with multiplication (a base point
in the operadic sense) admits an action of an operad weakly equivalent to the little disks operad.
As an immediate consequence, there exists a natural bracket on the (rational) homology of $\Tot \calO^{\bullet}$,
called the {\it Browder operation}.

On the other hand, Bousfield \cite{Bousfield87} constructed a spectral sequence computing the homology of a
totalization (under some conditions).
When moreover the cosimplicial space arises from the operad $\calO$ with multiplication,
then its $E^1$-term is the {\it Hochschild complex} of the homology operad $H_* (\calO )$.
It is known (Gerstenhaber--Voronov \cite{GerstenhaberVoronov95}, Turchin \cite{Tourtchine04,Tourtchine04_2})
that there exist a natural product and a bracket on such a complex which induce the {\it Gerstenhaber algebra structure},
the degree one Poisson algebra structure, on the homology.
Note that the Gerstenhaber algebra structure also comes from the action of
the chains of the little disks operad ({\it Deligne's conjecture}; see McClure--Smith \cite{McClureSmith02}).

The main result of this article states that the above two actions correspond with each other at least on the
homology level, or the two brackets coincide with each other.
Namely, for any operad $\calO$ with multiplication, Bousfield's spectral sequence computes $H_* (\Tot \calO^{\bullet} )$
as a Poisson algebra.
In some cases this gives us a method to compute the (topological) Browder operation purely algebraically,
as was done by Turchin \cite{Tourtchine04_2}.

In particular we concentrate on the space $\K$ of long knots in $\R^n$, $n>3$, or its variant $\K'$ (for definition
see Sinha \cite{Sinha04} or \fullref{2} below).
For $\K'$, Sinha \cite{Sinha04} constructed a model using the {\it Kontsevich operad}, essentially equivalent
to the little balls operad, based on the {\it embedding calculus} due to Goodwillie \cite{Goodwillie92,Goodwillie03}
and Goodwillie--Weiss \cite{GoodwillieWeiss99,Weiss99}.
In this case Bousfield's spectral sequence rationally degenerates at $E^2$ because of the {\it formality} of the
operad (see Kontsevich \cite{Kontsevich99}, Lambrechts--Turchin--Voli\'c \cite{LTV} and Lambrechts--Voli\'c \cite{LV}).
Thus the homology of the space of long knots is isomorphic to the Hochschild homology of the Kontsevich operad.
Our result shows that this is an isomorphism of Poisson algebras.

To prove the main result, we explicitly write down McClure--Smith
action, which compares with Budney's (possibly another) one
\cite{Budney03} defined on certain embedding spaces. In the case of
long knots, our explicit description suggests that McClure--Smith
action might be equivalent to Budney's one, but no rigorous proof of
their consistency has been given.  Budney's action is quite geometric;
``shrink one (framed) knot and make it go through another knot.''
Budney proved that the space of (framed) long knots in $\R^3$ is a
free 2-cubes object, hence the homology operations are highly
non-trivial (as studied by Budney and Cohen \cite{BudneyCohen05}).  In
higher dimensional case we might be able to approach the similar
freeness problem by means of the spectral sequence.

Note that Salvatore has announced our result \cite[Proposition
22]{Salvatore06}, but not written the proof. He proved that the space
$\K$ ($n>3$) is a double loop space but the projection $\K' \to \K$ is
not a double loop map, by use of the main result of this paper.  Here
we emphasize that our result implies the non-triviality of the
(topological) Browder operation, and this produces a homology class of
$\K$ which does not directly correspond to any chord diagram (see
\fullref{3}).

%%% 2nd section %%%

\section{The spaces of knots and the result}\label{2}

A {\it long knot} is an embedding $\R^1 \hookrightarrow \R^n$, $n\ge 3$, which agrees with a fixed line outside
a compact set. A long immersion is defined similarly.

Let $\K$ and $\Imm$ be the spaces of long knots and long immersions in $\R^n$ respectively. Consider the space
\[
 \K' :=\text{the homotopy fiber of }\K \hookrightarrow \Imm .
\]
Sinha's model for the space is as follows.

\begin{thm}{\rm \cite{Sinha04}}\qua
There exists a topological operad $X_n =\{ X_n (k)\}_{k\ge 0}$ with multiplication such that
\begin{enumerate}
\item there exist homotopy equivalences from $\mathrm{Conf}\, (\R^n ,k)=(\R^n )^k \setminus \Delta$ to $X_n (k)$
 for all $k\ge 0$,
\item when $n>3$, the homotopy invariant totalization $\tTot \X$ is weakly equivalent to $\K'$ where $\X$ is the
associated cosimplicial space.\qed
\end{enumerate}
\end{thm}

The operad $X_n$ is called the {\it Kontsevich operad} by Sinha \cite{Sinha04}.
Since it is multiplicative, we can construct an associated cosimplicial space (see \fullref{conventions}),
denoted by $\X$ in the above.

Below we deal with general cases. For any cosimplicial space $Y^{\bullet}$, its homotopy invariant totalization
is defined as
\[
 \tTot Y^{\bullet} = \Map (\tilde{\Delta}^{\bullet},Y^{\bullet} )
\]
where $\tilde{\Delta}^{\bullet}$ is a cofibrant replacement of the standard cosimplicial space $\Delta^{\bullet}$.
The ordinary totalization $\Tot Y^{\bullet}$ is defined by replacing $\tilde{\Delta}^{\bullet}$ by
$\Delta^{\bullet}$ in the above.
There is a map $\Tot Y^{\bullet} \to \tTot Y^{\bullet}$ induced by a canonical map
$\tilde{\Delta}^{\bullet} \to \Delta^{\bullet}$.

In the following let $\calO$ be an operad with multiplication, assuming $\calO (0)=\{ *\}$.
Denote by $\calO^{\bullet}$ its associated cosimplicial space.
One of the remarkable features of such a cosimplicial space $\calO^{\bullet}$ is the existence of a little disks
action on the totalization.

\begin{thm}{\rm \cite{McClureSmith04_2}}\qua
Let $\calO$ be a topological operad with multiplication, and $\calO^{\bullet}$ its associated cosimplicial space.
Then $\tTot \calO^{\bullet}$ admits an action of the operad $\calD$ which is weakly equivalent to the little disks operad.
That is, there exist maps
\[
 \gamma_k \co \calD (k) \times (\tTot \calO^{\bullet} )^k \longrightarrow \tTot \calO^{\bullet} ,\quad k\ge 1,
\]
satisfying the associativity conditions. Similar action exists on $\Tot \calO^{\bullet}$.\qed
\end{thm}

We denote the induced Browder operation by $\lambda$:
\[
 \lambda \co H_p (\tTot \calO^{\bullet} ,\Q )\otimes H_q (\tTot \calO^{\bullet} ,\Q )
 \longrightarrow H_{p+q+1}(\tTot \calO^{\bullet} ,\Q )
\]
We regard $\lambda$ as a Poisson algebra structure on the homology (see Cohen \cite{Cohen533}).

On the other hand, $H_* (\tTot \calO^{\bullet} )$ can be computed by means of Bousfield spectral sequence
\cite{Bousfield87}.
To do this, we choose a fibrant replacement $R:\calO^{\bullet}\we R\calO^{\bullet}$ of the
cosimplicial space $\calO^{\bullet}$.
$\Tot R\calO^{\bullet}$ may not be acted on by little disks since $R\calO^{\bullet}$ need not come from an operad.
But there exists a sequence of weak equivalences
\[
 \tTot \calO^{\bullet} \we\tTot R\calO^{\bullet} \wel \Tot R\calO^{\bullet} ,
\]
so there exists a Poisson algebra structure on $H_* (\Tot R\calO^{\bullet} )$ (we also call the bracket in this case
the Browder operation), induced from $H_* (\tTot \calO^{\bullet})$.

If $R\calO^{\bullet}$ satisfies the convergence conditions \cite[Theorem 3.2]{Bousfield87}, then the filtration
(defined below) on a double complex $S_* (R\calO^{\bullet} )$, the modules of singular chains of $R\calO^{\bullet}$,
yields a second quadratic spectral sequence $E^r$ converging to $H_* (\Tot R\calO^{\bullet} )$.
Indeed Bousfield proved that (under the same conditions) the total complex $TS(R\calO^{\bullet} )$ is an algebraic
model for $\Tot R\calO^{\bullet}$, that is, there exists a quasi-isomorphism
\[
 \varphi \co S_* (\Tot R\calO^{\bullet} ) \longrightarrow TS(R\calO^{\bullet} )
\]
(for its definition see \cite[section 2]{Bousfield87} or \fullref{4} below).
Its $E^1$-term is the Hochschild complex (see Turchin \cite{Tourtchine04,Tourtchine04_2}) of the homology operad
$H_* (R\calO^{\bullet} )$ (though $R\calO^{\bullet}$ may not be an operad, the natural isomorphism
$H_* (\calO^{\bullet} ) \xrightarrow{\cong} H_* (R\calO^{\bullet} )$ allows us to say that $H_* (R\calO^{\bullet})$
is an operad).
Since $\calO^{\bullet}$ is an operad with multiplication, there exist a natural product and a bracket
on $E^1$-term (see Gerstenhaber--Voronov \cite{GerstenhaberVoronov95}, Turchin \cite{Tourtchine04,Tourtchine04_2})
that make the $E^2$-term $HH_* (\calO^{\bullet})$ a Gerstenhaber algebra, that is, a Poisson algebra of degree one
\cite{Tourtchine04,Tourtchine04_2}.
This Poisson structure defined on $E^2$-term in fact descends on further terms (\fullref{spec_of_Poisson}).
Hence $E^{\infty}$-term also becomes a Poisson algebra.

The natural filtrations are defined by
\begin{gather*}
 F_p TS(R\calO^{\bullet} )|_{\deg =q} = \prod_{l\ge p}S_{q+l}(R\calO^l ), \\
 F_p H_q (TS(R\calO^{\bullet} )) = \img \left\{ H_q (F_p TS(R\calO^{\bullet} ))\to H_q (TS(R\calO^{\bullet} ))\right\} .
\end{gather*}
Denote its associated quotient by $GH_* (TS(R\X ))$:
\[
 GH_q (TS(R\calO^{\bullet} ))=\bigoplus_p F_p H_q (TS(R\calO^{\bullet} )) / F_{p+1}H_q (TS(R\calO^{\bullet} ))
\]
$H_* (\tTot \calO^{\bullet})$ is also filtered via isomorphism $\varphi^{-1}$. We will see
(\fullref{Browder_Tot}) that the Browder operation $\lambda$ preserves this filtration, hence via $\varphi$
a Poisson structure passes on $GH_* (TS(R\calO^{\bullet} ))$.

Thus we have two Poisson algebra structures on $GH_* (TS(R\calO^{\bullet} ))$: One induced via $\varphi$
and the other arising from the Hochschild homology. Our main result states that they are same.

\begin{thm}\label{main}
Let $\calO^{\bullet}$ be the cosimplicial space associated with a multiplicative operad $\calO$ satisfying the
convergence condition \cite[Theorem 3.2]{Bousfield87}. Choose a fibrant replacement $R\calO^{\bullet}$.
Then the Browder operation induced on $GH_* (TS(R\calO^{\bullet} ))$ via $\varphi$ coincides up to sign with
the Gerstenhaber structure defined on the Hochschild homology.
\end{thm}

In the case of $\calO =X_n$ we have the rational degeneracy of the spectral sequence.

\begin{thm}{\rm \cite{LTV,LV}}\qua
When $n>3$, the above homology spectral sequence degenerates at $E^2$-term rationally. Thus $H_* (\K' ,\Q )$ is
the Hochschild homology of the Poisson algebras operad.\qed
\end{thm}

This theorem holds since the {\it formality} of the operad $X_n$ allows us to replace the vertical differentials of
the double complex by zero maps. Hence we can compute the Browder operation on $H_* (\K' ,\Q )$ by calculating
$E^2$-term, without referring \fullref{spec_of_Poisson}.

The proof of \fullref{main} will be given later. First we compute some homology groups of $\K'$.

%%% 3rd section %%%

\section{Some computations}\label{3}

\begin{figure}[hbt]
\begin{center}
\includegraphics{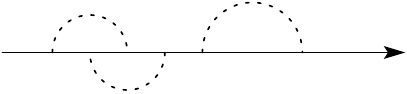}
\end{center}
\caption{an example of a chord diagram}
\label{chord_diagram}
\end{figure}

\begin{figure}[hbt]
\begin{center}
\includegraphics{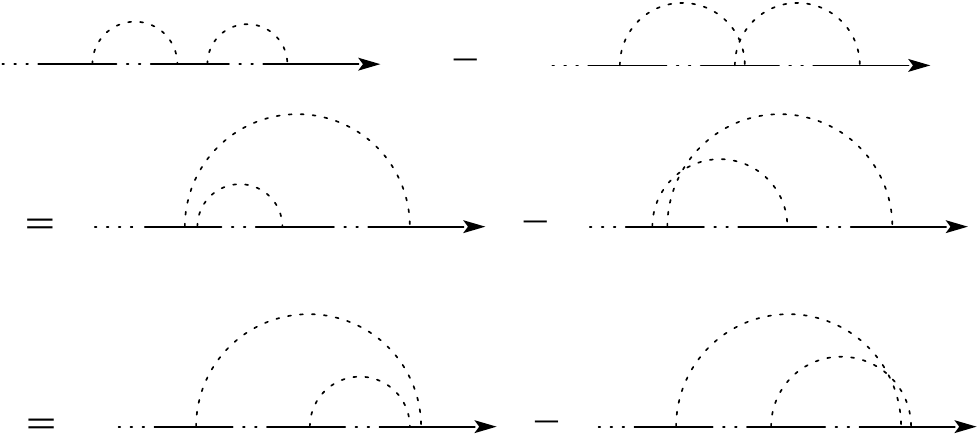}
\end{center}
\caption{4-term relations}
\label{4t}
\end{figure}

Below we assume $n>3$. The subgroup $\bigoplus_{k\ge 0}H_{(n-3)k}(\K )$ is known to be non-trivial
(Cattaneo--Cotta-Ramusino--Longoni \cite{CCL02}, Sinha \cite{Sinha04}, Turchin \cite{Tourtchine04_2}),
since it contains the subalgebra isomorphic to the algebra $\A$ of chord diagrams modulo 4-term and 1-term relations,
which we now explain.
An example of a chord diagram and the 4-term relator can be seen in Figures \ref{chord_diagram} and \ref{4t}
respectively. Since we consider the long knots, the chords are on the line, not on a circle.

By 1-term relation we regard a diagram $\Gamma$ as zero if $\Gamma$ has an {\it isolated chord}. Here a chord $c$
with endpoints $a,b\in \R^1$ ($a<b$) is said to be isolated if there is no other chord one of whose endpoints is in
$(a,b)$ and another outside $[a,b]$.

The degree of a chord diagram is $(n-3)k$ if it has $k$ chords.
The space of chord diagrams forms a graded algebra, whose product is defined as the concatenation of the diagrams.
Denote by $\A$ the algebra of chord diagrams modulo 4-term and 1-term relations.

\begin{figure}[hbt]
\begin{center}
\includegraphics[width=.7\hsize]{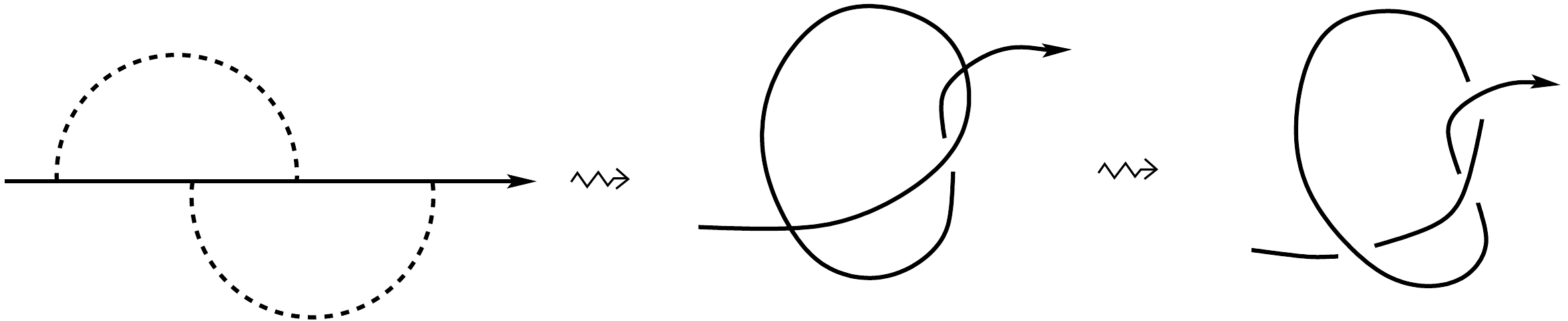}
\end{center}
\caption{the correspondence $\Gamma \mapsto \alpha (\Gamma )$}\label{blowup}
\end{figure}

Given a chord diagram $\Gamma$ with $k$ chords, we have a long immersion $f_{\Gamma}$ with $k$ transversal
self-intersections determined by the chords of $\Gamma$ (see \fullref{blowup}).
At each self-intersection we have resolutions of the intersection parametrized by $S^{n-3}$.
Considering all the resolutions of $k$ self-intersections, we have a map
\[
 s(\Gamma ) \co (S^{n-3})^k \longrightarrow \K .
\]
More explicitly, the knot $s(\Gamma )(z_1 ,\dots ,z_k )$ is defined as follows (see Cattaneo--Cotta-Ramusino--Longoni
\cite{CCL02});
if the $i$-th doublepoint of $f_{\Gamma}$ is $f_{\Gamma}(s_i )=f_{\Gamma}(u_i )$ ($s_i < u_i$), then
\[
 s(\Gamma )(z_1 ,\dots ,z_k ) (t)
 = f_{\Gamma}(t)+\delta z_i \exp \left[ \frac{1}{(t-s_i )^2 -\varepsilon^2}\right]
\]
for $t\in (s_i -\varepsilon ,\, s_i +\varepsilon )$, where $\delta$ and $\varepsilon$ are small positive numbers,
and otherwise $s(\Gamma )(z)(t) = f_{\Gamma}(t)$.

This map determines an $(n-3)k$-cycle $\alpha (\Gamma )$ of $\K$.
The correspondence $\alpha$ preserves the 4-term relations.

\begin{prop}{\rm \cite{CCL02,Tourtchine04_2}}\qua
The correspondence $\Gamma \mapsto \alpha (\Gamma )$ determines an injective homomorphism of algebras
\[
 \alpha \co \A \longrightarrow \bigoplus_{k\ge 0}H_{(n-3)k}(\K ). \rlap{\hspace{1.5in} \qedsymbol}
\]
\end{prop}

Thus we have a subalgebra of $H_* (\K )$ which is isomorphic to $\A$ and concentrated in the degrees $(n-3)k$, $k\ge 0$.
But few homology classes except those in $\alpha (\A )$ have been identified.

\fullref{main} implies that the spectral sequence computes the Poisson algebra structure on $H_* (\K' ,\Q )$.
By computing the $E^2$-term explicitly (see Turchin \cite{Tourtchine04_2}), we find that the Browder operation produces
new homology classes other than those in $\alpha (\A )$.

To state the result, we need the following.

\begin{lem}{\rm \cite{Sinha04}}\qua
For any $n\ge 3$, the inclusion $\K \hookrightarrow \Imm$ is null-homotopic. Consequently there exists a
homotopy equivalence $\K' \simeq \K \times \Omega^2 S^{n-1}$.\qed
\end{lem}

The topology of $\Omega^2 S^{n-1}$ is well known:
\[
 H_* (\Omega^2 S^{n-1},\Q )\cong
 \begin{cases}
  \Lambda^* [\lambda_{\Omega^2 S^{n-1}}(\iota ,\iota )]\otimes P[\iota ] & n\text{ is odd}, \\
  \Lambda^* [\iota ] & n\text{ is even},
 \end{cases}
\]
where $\iota$ is the image of $1\in \pi_{n-3}(\Omega^2 S^{n-1})$ via the Hurewicz isomorphism, and
$\lambda_{\Omega^2 S^{n-1}}$ is the Browder operation induced by the little disks action on $\Omega^2 S^{n-1}$.
Using this, we can determine the generators of some low-degree homology groups of $\K'$ and $\K$.

\begin{figure}[hbt]
\begin{center}
\includegraphics[scale=.9]{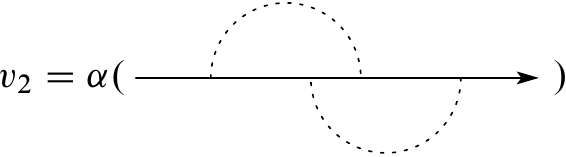}
\end{center}\vspace{-.5cm}
\caption{the element $v_2$}
\label{v_2}
\end{figure}

\begin{thm}
Suppose $n>4$. Then
\[
 H_{3n-8}(\K' ,\Q ) \cong
  \begin{cases}
   \Q^2 & n\text{ is odd}, \\
   \Q   & n\text{ is even}.
  \end{cases}
\]
When $n$ is odd, the above group is generated by $\iota \lambda (\iota ,\iota )$ and $\lambda (\iota ,v_2 )$,
where $v_2 \in H_{2(n-3)}(\K )$ is the class made from a chord diagram with two interleaving chords
(see \fullref{v_2}). By the K\"unneth formula
\[
 H_{3n-8}(\K ,\Q )\cong \Q
\]
and this group is generated by the homology class corresponding to $\lambda (\iota ,v_2 )$, which is the first
example of a homology class of $\K$ which does not directly correspond to chord diagrams.\qed
\end{thm}

In other words, we can regard $H_* (\Omega^2 S^{n-1},\Q )$ as acting on $H_* (\K ,\Q )$ via the Browder operation,
and this action is non-trivial when $n$ is odd. The homology class $\lambda (\iota ,v_2 )$ (when $n$ is odd) is not
in $\alpha (\A )$, since it is not in degrees $(n-3)k$.

When $n$ is even, $\iota$ belongs to the center (Turchin \cite{Tourtchine04_2}).
This is because the map of spectral sequences $E^2 (\K' )\to E^2 (\K )$ is a map of Gerstenhaber algebras, and
under the map $\iota$ must be sent to 0 for the dimensional reason.
But $\rank H_{3n-8}(\K' )=1$ still holds, and the generator is unknown yet.

Salvatore \cite{Salvatore06} pointed out that $\K$ is also a double loop space, but the existence of
$\lambda (\iota ,v_2 )$ shows that the projection $\K' \to \K$ does not preserve the Browder operation
(since $\iota$ is mapped to zero by the dimensional reason), hence this projection is not the map of double loop
spaces.

The above computation suggests that we might be able to obtain more homology classes by applying
$\lambda (\iota ,\cdot )$ to the elements of $\alpha (\A )$. It is possible in principle to compute them, but
it would become more exhausting as the degrees increase (see Cattaneo--Cotta-Ramusino--Longoni \cite{CCL02}, Turchin
\cite{Tourtchine04_2}).

%%% 4th section %%%

\section{Proof of the main theorem}\label{4}

\subsection{Disks action on $\Tot \calO^{\bullet}$}

First we let $\calO$ be a topological operad with multiplication and suppose the associated cosimplicial space
$\calO^{\bullet}$ is fibrant.

To obtain the explicit formula for the Browder operation on $H_* (\Tot \calO^{\bullet})$, it suffices to study
in detail the map
\[
 \gamma \co S^1 \times (\Tot \calO^{\bullet} )^2 \longrightarrow \Tot \calO^{\bullet}
\]
($S^1 \sim \calD (2)$) defined by McClure and Smith \cite{McClureSmith04_2}.
We think $g^{\bullet} \in \Tot \calO^{\bullet}$ as a sequence of maps $g^l :\Delta^l \to \calO^l$
($l \ge 0$), compatible with the cosimplicial structure maps. Defining the above map $\gamma$ is equivalent to
defining the maps
\[
 \gamma (\tau ; g^{\bullet}_1 ,g^{\bullet}_2 )^l \co \Delta^l \longrightarrow \calO^l, \quad l \ge 0
\]
compatible with cosimplicial structure maps, for any $\tau \in S^1$, $g^{\bullet}_i \in \Tot \calO^{\bullet}$.
After the fashion of McClure and Smith's work \cite{McClureSmith02}, we first look at the easiest example.

\subsubsection{Conventions}\label{conventions}

We now proceed to operadic computations.  For such notions see
McClure--Smith \cite{McClureSmith02,McClureSmith04_2}.  We use the
notation $\circ_i$ for operadic structure:
\[
 \circ_i \co \calO (p) \times \calO (q) \longrightarrow \calO (p+q-1),\quad
 a\circ_i b =a(\overbrace{\id ,\dots ,\id}^{i-1},b,\id ,\dots ,\id )
\]
for any operad $\calO$. The symbol $\circ_i$ represents the `inserting' operation.

Recall that an operad $\calO$ is said to be {\it multiplicative} if we can choose basepoints
$\mu_k \in \calO (k)$, $k \ge 0$, in the operadic sense;
\[
 \mu_k (\mu_{j_1} ,\dots ,\mu_{j_k}) = \mu_{j_1 +\dots +j_k}.
\]
In this case $\calO$ can be seen as a cosimplicial space by letting $\calO^k := \calO (k)$ and defining
\begin{alignat*}{2}
 d^i \co \calO^k &\longrightarrow \calO^{k+1}, & \quad & 0 \le i \le k+1, \\
 s^i \co \calO^k &\longrightarrow \calO^{k-1}, & & 0 \le i \le k-1,
\end{alignat*}
\begin{align*}
 d^i (x) &=\tag*{\hbox{by}}
 \begin{cases}
  \mu_2 \circ_2 x & i=0, \\
  x \circ_i \mu_2 & 1 \le i \le k, \\
  \mu_2 \circ_1 x & i=k+1,
 \end{cases} \\
 s^i (x) &= x \circ_{i+1} e
\end{align*}
where $e \in \calO (0)$. They indeed satisfy the cosimplicial identities.

Below by convention
\begin{align*}
 \Delta^k &:=
 \begin{cases}
  \{ *\} & k=0, \\
  \{ (t_1 ,\dots ,t_k ) \in [-1,1]^k \, |\, t_1 \le \dots \le t_k \} & k \ge 1,
 \end{cases} \\
 S^1 &:= [-1,1] / \sim .
\end{align*}

\subsubsection{The first case}\label{l=0}

\begin{figure}[htb]
\begin{center}
\includegraphics{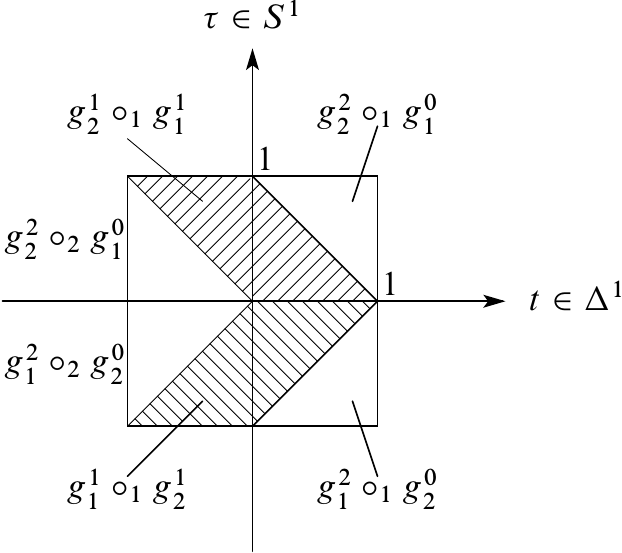}
\end{center}\vspace*{-.5cm}
\caption{the `definition' of $\gamma (\tau ; g^{\bullet}_1 ,g^{\bullet}_2 )^1 (t)$}\label{gamma}
\end{figure}

As the easiest case we construct
\[
 \gamma (\tau ; g^{\bullet}_1 ,g^{\bullet}_2 )^1 \co \Delta^1 \longrightarrow \calO^1
\]
(assuming $\calO^0 =\{ e\}$, the case $l=0$ is obvious).
Imitating a work of McClure and Smith \cite{McClureSmith02}, we roughly illustrate the definition of this map
(see \fullref{gamma}).

\begin{rem}
Broadly speaking, as $\tau$ increases from $-1$ to $0$, the knot $g_2$ `goes through' the knot $g_1$
($\circ_i$ will represent the `insertion').
By $\tau =0$, $g_2$ will get away from $g_1$ and `juxtapose' to $g_1$
($\mu_2 \in \calO (2)$ will indicate `concatenation').
When $0< \tau <1$, now $g_1$ passes through $g_2$, and juxtaposes to $g_2$ when $\tau=1$.
See also the pictures in Budney's paper \cite{Budney03}.
The reason why we may regard $\circ_i$ and $\mu_2$ as respectively `insertion' and `concatenation' can be found
in the definition of {\it Poisson algebras operad} $H_* (X_n )$ (see Turchin \cite{Tourtchine04,Tourtchine04_2}).\qed
\end{rem}

More precise definition is as follows. When $-1 < \tau < 0$, define non-negative integers $i_0 ,i_1$ by
\[
 (i_0 ,i_1 )=
 \begin{cases}
  (1,1) & t<\tau ,\\
  (0,1) & \tau \le t <1+\tau , \\
  (0,0) & 1+\tau \le t,
 \end{cases}
\]
and define $(u_1 ,u_2 )\in \Delta^{2+i_0 -i_1} \times \Delta^{i_1 -i_0}$ by
\[
 (u_1 ,u_2 ) =
 \begin{cases}
  ((2t+1,1+2\tau );\, *) \in \Delta^2 \times \Delta^0 & t<\tau , \\
  (2\tau +1;\, 2(t-\tau )-1) \in \Delta^1 \times \Delta^1 & \tau \le t < 1+\tau , \\
  ((1+2\tau ,2t-1);\, *) \in \Delta^2 \times \Delta^0 & 1+\tau \le t.
 \end{cases}
\]
Then
\[
 \gamma (\tau ;g^{\bullet}_1 ,g^{\bullet}_2 )^1 (t) := g^{2+i_0 -i_1}_1 (u_1 ) \circ_{i_0 +1} g^{i_1 -i_0}_2 (u_2 )
 \in \calO^1 .
\]
When $0< \tau <1$, define
\[
 (j_0 ,j_1 ) =
 \begin{cases}
  (1,1) & t<-\tau , \\
  (0,1) & -\tau \le t < 1-\tau , \\
  (0,0) & 1-\tau \le t,
 \end{cases}
\]
and define $(v_1 ,v_2 )\in \Delta^{j_1 -j_0} \times \Delta^{2+j_0 +j_1}$ by
\[
 (v_1 ,v_2 ) =
 \begin{cases}
  (*;\, (2t+1,1-2\tau )) \in \Delta^0 \times \Delta^2 & t < -\tau , \\
  (2(t-\tau )-1;\, 1-2\tau ) \in \Delta^1 \times \Delta^1 & -\tau \le t < 1-\tau , \\
  (*;\, (1-2\tau ,2t-1)) \in \Delta^0 \times \Delta^2 & 1-\tau \le t.
 \end{cases}
\]
Then
\[
 \gamma (\tau ;g^{\bullet}_1 ,g^{\bullet}_2 )^1 (t) := g^{j_1 -j_0}_2 (v_2 ) \circ_{j_0 +1} g^{2+j_0 -j_1}_1 (v_1 ).
\]
The cases $\tau =0,\pm 1$ may be contained in the above definitions; even in those cases
$\gamma (\tau ;g^{\bullet}_1 ,g^{\bullet}_2 )^1 (t)$ is well-defined.
But it would be better to give the definitions for $\tau =0,\pm 1$ separately to make the meaning of $\gamma$
clearer.
We define
\begin{align*}
 \gamma (0 ;g^{\bullet}_1 ,g^{\bullet}_2 )^1 (t) &:= \mu_2 (g^{i_0}_1 (w_1),g^{1-i_0}_2 (w_2 )), \\
 \gamma (\pm 1 ;g^{\bullet}_1 ,g^{\bullet}_2 )^1 (t) &:= \mu_2 (g^{j_1}_2 (w_2 ),g^{1-j_1}_1 (w_1 )),
\end{align*}
where $i_0$ (for $\tau =0$) and $j_1$ (for $\tau =+1$) are as above, and
\[
 (w_1 ,w_2 ) =
 \begin{cases}
  (2t+1;\, *) \in \Delta^1 \times \Delta^0 & t<0, \\
  (*;\, 2t-1) \in \Delta^0 \times \Delta^1 & t \ge 0.
 \end{cases}
\]
Then $\gamma (\tau ;g^{\bullet}_1 ,g^{\bullet}_2 )^1 (t)$ is indeed well-defined.
For example, let us see the continuity at $t=1+\tau$ when $-1 < \tau < 0$. By definition
\[
 \gamma (\tau ;g^{\bullet}_1 ,g^{\bullet}_2 )^1 (t) =
 \begin{cases}
  g^1_1 (2\tau +1) \circ_1 g^1_2 (2(t-\tau )-1) & t< 1+\tau , \\
  g^2_1 (2\tau +1,2\tau +1) \circ_1 g^0_2 (*) & t= 1+\tau .
 \end{cases}
\]
Then we can show
\[
 \lim_{t \nearrow \tau +1} \gamma (\tau ;g^{\bullet}_1 ,g^{\bullet}_2 )^1 (t) =g^1_1 (2\tau +1) \circ_1 g^1_2 (1)
 =\gamma (\tau ;g^{\bullet}_1 ,g^{\bullet}_2 )^1 (1+\tau )
\]
as follows;
\begin{alignat*}{2}
 &g^1_1 (2\tau +1) \circ_1 g^1_2 (1) & \quad & \ \\
 &\quad = g^1_1 (2\tau +1) \circ_1 g^1_2 (d^1 (*)) = g^1_1 (2\tau +1) \circ_1 d^1 g^0_2 (*)
  & & (g^l d^i =d^i g^{l-1}) \\
 &\quad = g^1_1 (2\tau +1) (\mu_2 (g^0_2 (*),\id )) & & \\
 &\quad = g^1_1 (2\tau +1)(\mu_2 ) (g^0_2 (*), \id )
  & & \text{(associativity of } \calO ) \\
 &\quad = (d^1 g^1_1 (2\tau +1))(g^0_2 (*),\id ) = g^2_1 (d^1 (2\tau +1))(g^0_2 (*),\id ) & & \\
 &\quad = g^2_1 (2\tau +1, 2\tau +1) \circ_1 g^0_2 (*) = \gamma (\tau ;g^{\bullet}_1 ,g^{\bullet}_2 )^1 (1+\tau ).
  & &
\end{alignat*}

We see one more point; $\gamma$ is continuous at $\tau =0, \pm 1$. For example we show
\[
 \lim_{\tau \searrow -1}\gamma (\tau ;g^{\bullet}_1 ,g^{\bullet}_2 )^1 (t)
 = \gamma (-1;g^{\bullet}_1 ,g^{\bullet}_2 )^1 (t)
\]
when $-1<t<0$. In this case the above limit is equal to
\begin{align*}
 g^1_1 (-1) \circ_1 g^1_2 (2t+1) &= g^1_1 (d^0 (*))(g^1_2 (2t+1)) = (d^0 (g^0_1 (*)))(g^1_2 (2t+1)) \\
 &= \mu_2 (\id ,g^0_1 (*))(g^1_2 (2t+1)) \\
 &= \mu_2 (g^1_2 (2t+1),g^0_1 (*)) = \gamma (-1;g^{\bullet}_1 ,g^{\bullet}_2 )^1 (t).
\end{align*}
The proofs for other cases go in similar ways.

\subsubsection{General cases}

Here we describe the complete definition of the map
\[
 \gamma \co S^1 \times (\Tot \calO^{\bullet})^2 \longrightarrow \Tot \calO^{\bullet},
\]
a special case of the construction by McClure and Smith \cite{McClureSmith04_2}, and the induced map on homology.

First we define the integers $i_{\varepsilon},j_{\varepsilon}$ ($\varepsilon =0,1$), which determine where
$g_{\sigma (2)}$ is `inserted' in $g_{\sigma (1)}$, for any $t=(t_1 ,\dots ,t_l )\in \Delta^l$ and $\tau \in S^1$.
When $-1< \tau <0$, define
\[
 i_0 :=\text{min}\{ i;\, t_{i+1} \ge \tau \},\quad
 i_1 :=\text{min}\{ i;\, t_{i+1} \ge 1+\tau \},
\]
and when $0< \tau <1$, define
\[
 j_0 :=\text{min}\{ j;\, t_{j+1} \ge -\tau \} ,\quad
 j_1 :=\text{min}\{ j;\, t_{j+1} \ge 1-\tau \} .
\]

\begin{figure}[htb]
\begin{center}
\includegraphics{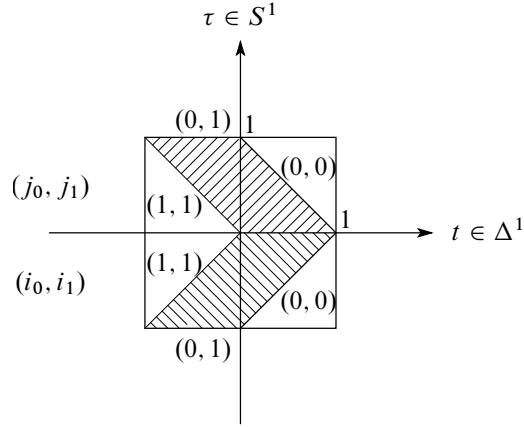}
\end{center}\vspace*{-.5cm}
\caption{the values $i_{\varepsilon},j_{\varepsilon}$ for $l=1$}\label{01}
\end{figure}

Next, define
\[
 (u_1 (\tau ,t),u_2 (\tau ,t))\in \Delta^{l+i_0 -i_1 +1}\times \Delta^{i_1 -i_0}
\]
(when $-1< \tau <0$) and
\[
 (v_1 (\tau ,t),v_2 (\tau ,t))\in \Delta^{j_1 -j_0}\times \Delta^{l+j_0 -j_1 +1}
\]
(when $0< \tau <1$) by
\begin{align*}
 u_1 (\tau ,t) &= (2t_1 +1 ,\dots ,2t_{i_0}+1, 1+2\tau , 2t_{i_1 +1}-1 ,\dots ,2t_l -1), \\
 u_2 (\tau ,t) &= (2t_{i_0 +1}-2\tau -1 ,\dots ,2t_{i_1}-2\tau -1), \\
 v_1 (\tau ,t)&= (2t_{j_0 +1}-2\tau -1 ,\dots ,2t_{j_1}-2\tau -1), \\
 v_2 (\tau ,t)&= (2t_1 +1 ,\dots ,2t_{j_0}+1, 1-2\tau , 2t_{j_1 +1}-1 ,\dots ,2t_l -1).
\end{align*}
Here we note a consequence of a straightforward computation.

\begin{lem}\label{coord}
For any given $l$ and $0\le i_0 \le i_1 \le l$, define
\[
 \Delta (i_0 ,i_1 ):= \{ (\tau ,t)\in (-1,0) \times \Delta^l \, |\,
 i_{\epsilon} (\tau ,t) =i_{\epsilon}\ (\epsilon =0,1) \} .
\]
Then the correspondence
\[
 u=(u_1 ,u_2 ) \co \Delta (i_0 ,i_1 )\longrightarrow \Delta^{l+i_0 -i_1 +1}\times \Delta^{i_1 -i_0},
\]
defined in the above remark, is a homeomorphism on their interior. Similarly, if we define
\[
 \Delta' (j_0 ,j_1 ):= \{ (\tau ,t)\in (0,1) \times \Delta^l \, |\,
 j_{\epsilon} (\tau ,t) =j_{\epsilon}\ (\epsilon =0,1) \} ,
\]
for any given $l$ and $0\le j_0 \le j_1 \le l$, then
\[
 v=(v_1 ,v_2 ) \co \Delta' (j_0 ,j_1 )\longrightarrow  \Delta^{j_1 -j_0}\times \Delta^{l+j_0 -j_1 +1}
\]
is a homeomorphism on their interior.\qed
\end{lem}

Using $u_i$ and $v_i$, we define
\[
 \gamma (\tau ;g^{\bullet}_1 ,g^{\bullet}_2 )^l (t)=
 \begin{cases}
  g^{l+i_0 -i_1 +1}_1 (u_1 (\tau ,t)) \circ_{i_0 +1} g^{i_1 -i_0}_2 (u_2 (\tau ,t)) & -1 < \tau < 0, \\
  g^{j_1 -j_0}_2 (v_2 (\tau ,t)) \circ_{j_0 +1} g^{l+j_0 -j_1 +1}_1 (v_1 (\tau ,t)) & 0 < \tau < 1.
 \end{cases}
\]

For $\tau =0,\pm 1$ (at which $g_1$ and $g_2$ `switch'), we separately give the definition; when $\tau =0$,
define $i_0$ as above and
\begin{align*}
 \bar{u}_1 (\tau ,t) &:= (2t_1 +1,\dots ,2t_{i_0}+1) \in \Delta^{i_0}, \\
 \bar{u}_2 (\tau ,t) &:= (2_{i_0 +1}-1,\dots ,2t_l -1) \in \Delta^{l-i_0}.
\end{align*}
Then we define
\[
 \gamma (0; g^{\bullet}_1 ,g^{\bullet}_2 )^l (t) := \mu_2 (g^{i_0}_1 (\bar{u}_1 ),g^{l-i_0}_2 (\bar{u}_2)).
\]
When $\tau =\pm 1$, define $j_1$ as above (for $\tau = +1$) and
\begin{align*}
 \bar{v}_1 (\tau ,t) &:= (2t_{j_1 +1}-1,\dots ,2t_l -1) \in \Delta^{l-j_1}, \\
 \bar{v}_2 (\tau ,t) &:= (2t_1 +1 ,\dots ,2t_{j_1}+1) \in \Delta^{j_1}.
\end{align*}
Then we define
\[
 \gamma (\pm 1 ;g^{\bullet}_1 ,g^{\bullet}_2 )^l (t) := \mu_2 (g^{j_1}_2 (\bar{v}_2 ), g^{l-j_1}_1 (\bar{v}_1 )).
\]
Then the map $\gamma (\tau ;g^{\bullet}_1 ,g^{\bullet}_2 )^l (t)$ is well-defined, which is proven similarly as when
$l=1$.
Moreover $\gamma$ indeed defines a map to a totalization.

\begin{lem}
The sequence of maps $\gamma (\tau ;g^{\bullet}_1 ,g^{\bullet}_2 )^{\bullet}$ is compatible with the cosimplicial
structure maps.\qed
\end{lem}

For example, if $(\tau ,t)\in \Delta^l (i_0 ,i_1 )$, $0 < i_0 < i_1$ and $0 < i \le i_0$, then
$d^i (t)=(\dots ,t_i ,t_i ,\dots ,t_{i_0},\dots )$ and $(\tau ,d^i (t))\in \Delta^{l+1}(i_0 +1,i_1 +1)$. Hence
\begin{align*}
 &\gamma (\tau ; g^{\bullet}_1 ,g^{\bullet}_2 )^{l+1}(d^i (t)) \\
 &\quad =g_1 (\dots ,2t_i +1 ,2t_i +1,\dots ) \circ_{i_0 +2} g_2 (2t_{i_0 +1}-2\tau -1,\dots ,2t_{i_1}-2\tau -1) \\
 &\quad =g_1 (d^i (\dots ,2t_i +1 ,\dots )) \circ_{i_0 +2} g_2 (u_2 (\tau ,t)) \\
 &\quad = (d^i g_1 (u_1 (\tau ,t))) \circ_{i_0 +2}g_2 (u_2 (\tau ,t))
  = (g_1 (u_1 ) \circ_i \mu_2 )\circ_{i_0 +2} g_2 (u_2 ) \\
 &\quad = (g_1 (u_1 ) \circ_{i_0 +1} g_2 (u_2 )) \circ_i \mu_2
  = d^i (\gamma (\tau ; g^{\bullet}_1 ,g^{\bullet}_2 )^l (t))
\end{align*}
(the fifth equality uses $i \le i_0$). Other cases can be proven similarly.

\begin{rem}
In fact the definitions for $\tau =0,\pm 1$ can be obtained as the limits of those for $-1 < \tau <0$ and
$0 < \tau <1$. But we give them separately to clarify the meaning of the definition.\qed
\end{rem}

For the induced map on homology, we only need to pre-compose the Eilenberg-MacLane map (see \fullref{alg_model}).
In the following we use the same symbols as above.

\begin{thm}\label{action}
Let $g^{\bullet}_1 \in S_q (\Tot \calO^{\bullet})$ and $g^{\bullet}_2 \in S_s (\Tot \calO^{\bullet})$ be cycles.
Define the map
\[
 \ang{g^{\bullet}_1}{g^{\bullet}_2}^l \co \Delta^q \times \Delta^s \times \Delta^1 \times \Delta^l
 \longrightarrow \calO^l ,\quad l\ge 0
\]
by
\[
 \ang{g^{\bullet}_1}{g^{\bullet}_2}^l (x,y,\tau ,t) :=
 \begin{cases}
  \mu_2 \left( g^{j_1}_2 (y)(\bar{v}_2 ),\ g^{l-j_1}_1 (x)(\bar{v}_1 ) \right) & \tau =\pm 1, \\
  g^{l+i_0 -i_1 +1}_1 (x)(u_1 )\circ_{i_0 +1}g^{i_1 -i_0}_2 (y)(u_2 ), & -1 < \tau < 0, \\
  \mu_2 \left( g^{i_0}_1 (x)(\bar{u}_1 ),\ g^{l-i_0}_2 (y)(\bar{u}_2 ) \right) & \tau =0, \\
  g^{l+j_0 -j_1 +1}_2 (y)(v_2 )\circ_{j_0 +1} g^{j_1 -j_0}_1 (x)(v_1 ), & 0 < \tau < 1.
 \end{cases}
\]
Then the map
\[
 \lambda \co H_q (\Tot \calO^{\bullet})\otimes H_s (\Tot \calO^{\bullet})\to H_{q+s+1}(\Tot \calO^{\bullet} )
\]
given by
\[
 \lambda (g^{\bullet}_1 ,g^{\bullet}_2 ):=(-1)^{q+1}\ang{g^{\bullet}_1}{g^{\bullet}_2}\circ EM
 \in H_{q+s+1}(\Tot \calO^{\bullet} )
\]
is the Browder operation, where $EM$ is the Eilenberg-MacLane map.\qed
\end{thm}

\subsection{Algebraic model for $\Tot \calO^{\bullet}$}\label{alg_model}

Let $TS(\calO^{\bullet})$ be the total complex of the double complex from \S 2: its degree $k$ part is
\[
 TS(\calO^{\bullet} )_k =\prod_{l\ge 0}S_{k+l}(\calO^l ),
\]
and the differential is
\[
 \partial_T =\partial +(-1)^p \delta,
\]
where $\partial =\sum (-1)^i d^i_*$ ($d^i$ was defined in \fullref{conventions}), a signed sum of the coface maps,
and $\delta$ is the usual boundary map of singular chain complex.
Bousfield \cite{Bousfield87} constructed a quasi-isomorphism
\[
 \varphi \co S_* (\Tot \calO^{\bullet} )\longrightarrow TS(\calO^{\bullet} )
\]
when $\calO^{\bullet}$ satisfies some conditions. This is defined as follows. We regard a chain
$f^{\bullet}=\sum a_i f^{\bullet}_i \in S_q (\Tot \calO^{\bullet} )$ as the sum of maps
$f^{\bullet}_i =\{ f^l_i \}_{l\ge 0}$,
\[
 f^l_i \co \Delta^q \times \Delta^l \longrightarrow \calO^l ,
\]
which is compatible with the cosimplicial structure maps of the $\Delta^l$-factor. We choose the Eilenberg-MacLane map
$EM\in S_{q+l}(\Delta^q \times \Delta^l )$, $l,q\ge 0$, which gives a chain equivalence
\[
 S_* (M)\otimes S_* (N)\xrightarrow{\simeq} S_* (M\times N)
\]
for any spaces $M$ and $N$. Then the quasi-isomorphism $\varphi$ is defined by
\begin{gather*}
 \varphi(f^{\bullet})\in \prod_{l\ge 0}S_{q+l}(\calO^l ), \\
 \varphi (f^{\bullet})_l :=\sum a_i (f^l_i \circ EM)\in S_{q+l}(\calO^l ).
\end{gather*}
where we write $h=(h_l )_{l\ge 0}$ with $h_l \in S_{q+l}(\calO^l )$ for any $h\in TS(\calO^{\bullet})_q$.
Our main theorem states that the induced isomorphism on homology preserves the Poisson algebra structures.

The spectral sequence associated with the filtration defined in \fullref{2} on the double complex converges
strongly when $E^1_{-p,q}=0$ for $p>q$ and, for any $m \ge 0$, there are only finitely many $(p,q)$ such that
$q-p=m$ and $E^1_{-p,q}\ne 0$. In the case of Kontsevich operad $X_n$, it turns out that the `normalized' $E^1$-term
\[
 E^1 \cap \bigcap \ker s^i_*
\]
satisfies those conditions (for definition of $s^i$ see \fullref{conventions}).
The proof uses the explicit form of the homology of $X_n (k) \simeq \mathrm{Conf}\, (\R^n ,k)$ (see Sinha \cite{Sinha04}).

\subsection{Browder operation in terms of $TS(\calO^{\bullet} )$}

Via the quasi-isomorphism $\varphi$, the Browder operation will be interpreted as follows.

\begin{thm}\label{Browder_Tot}
For any $g^{\bullet}_1 \in H_q (\Tot \calO^{\bullet})$ and $g^{\bullet}_2 \in H_s (\Tot \calO^{\bullet})$,
\begin{multline*}
 \varphi (\lambda (g^{\bullet}_1 ,g^{\bullet}_2 ))_l = \sum_{p+r=l+1}(-1)^{(p+1)(r+1)+q+s}
 \left[ \sum_{i=1}^p (-1)^{\epsilon_i} \varphi (g^{\bullet}_1 )_p \circ_i \varphi (g^{\bullet}_2 )_r \right. \\
 -\left. (-1)^{(q+1)(s+1)}
 \sum_{j=1}^r (-1)^{\epsilon'_j}\varphi (g^{\bullet}_2 )_r \circ_j \varphi (g^{\bullet}_1 )_p \right]
\end{multline*}
where
\[
 \epsilon_i =(q-1)(p-i)+(p-1)(r+s),\quad \epsilon'_j =(p-1)(r-j)+(p+q)(r-1).
\]
\end{thm}

\begin{proof}
The map
\[
 \Delta^q \times \Delta^s \times \Delta^1 \times \Delta^l \xrightarrow{\ang{g_1}{g_2}^l} X^l
\]
pre-composed by $EM$ represents $\varphi (\lambda (g_1 ,g_2 ))_l$. By definition $\Delta^1 \times \Delta^l$ is
decomposed by $\Delta (i_0 ,i_1 )$'s and $\Delta' (j_0 ,j_1 )$'s ($0\le i_0 \le i_1 \le l$, $0\le j_0 \le j_1 \le l$),
see \fullref{01}.
By \fullref{action}, $\ang{g_1}{g_2}^l$ is $g^p_1 \circ_{i_0 +1}g^r_2$ when restricted on
$\Delta^q \times \Delta^s \times \Delta (i_0 ,i_1 )$, where
\[
 p=l+i_0 -i_1 +1,\quad r=i_1-i_0 ,
\]
and, when restricted on $\Delta^q \times \Delta^s \times \Delta' (j_0 ,j_1 )$,
$\ang{g_1}{g_2}^l$ is $g^r_2 \circ_{j_0 +1}g^p_1$, where $p$ and $r$ are determined similarly.

Thus $\varphi (\lambda (g^{\bullet}_1 ,g^{\bullet}_2 ))$ should be a linear sum of
$\varphi (g^{\bullet}_1 )_p \circ_i \varphi (g^{\bullet}_2 )_r$, $1\le i\le p$, and
$\varphi (g^{\bullet}_2 )_r \circ_j \varphi (g^{\bullet}_1 )_p$, $1\le j\le r$. The coefficients are $\pm 1$
because of \fullref{coord}. The signs $\pm$ are those of the Jacobians of the maps
\begin{gather*}
 \Delta^q \times \Delta^s \times \Delta (i_0, i_1 ) \xrightarrow{\approx}
 (\Delta^q \times \Delta^p )\times (\Delta^s \times \Delta^r ), \\
 \Delta^q \times \Delta^s \times \Delta' (j_0, j_1 ) \xrightarrow{\approx}
 (\Delta^s \times \Delta^r )\times (\Delta^q \times \Delta^p ),
\end{gather*}
given by $u,v$ as in \fullref{action}. Explicit computations show that the signs are
\begin{gather*}
 (-1)^{1+i_0 (i_1 -i_0 )(l-i_1 )+ps}=(-1)^{(p+1)(r+1)+s+1}(-1)^{(r-1)(p-i_0 -1)+(p-1){r+s}},\\
 (-1)^{1+j_0 +(j_1 -j_0 )(l-j_1 )+q(r+s)}=(-1)^{(p+1)(r+1)+qs+q+1}\times
\tag*{\rm and}\\
\hspace{7cm}(-1)^{(p-1)(r-j_0 -1)+(p+q)(r-1)}
\end{gather*}
respectively. They give the desired formula.
\end{proof}

The isomorphism $\varphi$ introduces a filtration on $H_* (\Tot \calO^{\bullet})$;
\[
 F_p H_* (\Tot \calO^{\bullet}) := \varphi^{-1}F_pH_* (TS(\calO^{\bullet})).
\]
The definition of $F_*$ together with \fullref{Browder_Tot} says that the Browder operation $\lambda$
preserves this filtration in the sense
\[
 x\in F_p ,\ y \in F_r \ \Longrightarrow \ \lambda (x,y) \in F_{p+r-1}.
\]

\subsection{Poisson bracket on $E^2$}

Here let $\calO'$ be any operad of graded modules with multiplication (in our case $\calO'$ will be $H_* (\calO )$).
By unraveling the descriptions of the Hochschild complex $(\calO' ,\partial )$ (see Gerstenhaber--Voronov
\cite{GerstenhaberVoronov95}, Turchin \cite{Tourtchine04,Tourtchine04_2}), we can see the following.

\begin{thm}\label{turchin}{\rm \cite{GerstenhaberVoronov95,Tourtchine04,Tourtchine04_2}}\qua
The Poisson structure on $H_* (\calO' )$ is induced by the (degree-preserving) map
\begin{gather*}
 \Psi \co \calO' (p)\times \calO' (r)\longrightarrow \calO' (p+r-1), \\
 \Psi (x,y) := \sum_{1\le i\le p}(-1)^{\epsilon_i}x\circ_i y
 -(-1)^{(q+1)(s+1)}\sum_{1\le j\le r}(-1)^{\epsilon'_j}y\circ_j x
\end{gather*}
where $q=\deg x-p$, $s=\deg y-r$, and $\epsilon_i$, $\epsilon'_j$ are as in \fullref{Browder_Tot}.\qed
\end{thm}

In our case $\Psi$ is defined on $E^1$ and makes $E^2$ a Poisson algebra.
Indeed $E^r$ is a spectral sequence of a Poisson algebra because of the following.

\begin{prop}\label{spec_of_Poisson}
We have $d^r \Psi (x,y)=\Psi (d^r x,y)+(-1)^{\abs{x}+r}\Psi (x,d^r y)$ on $E^r$, $r\ge 2$.\qed
\end{prop}

The proof uses the definition of the boundary operation of the Hochschild complex (see \fullref{alg_model}).

Thus $E^{\infty}$ inherits the induced Poisson bracket, and via the isomorphism
\[
 \psi :E^{\infty}_{-p,q} \longrightarrow (F_p /F_{p+1})H_q (TS(\calO^{\bullet})),
\]
$GH_* (TS(\calO^{\bullet}))$ also becomes a Poisson algebra, where $G$ denotes the associated quotient.
Comparing \fullref{Browder_Tot} with \fullref{turchin}, we can see that this bracket coincides
on $GH_* (TS(\calO^{\bullet}))$ with the Browder operation induced via $\varphi$.

\subsection{The case of the space of knots}

Finally consider the case that $\calO^{\bullet}$ is not fibrant (in particular the case of the space of knots,
$\calO =X_n$).
Though we must use the fibrant replacement $R\calO^{\bullet}$, this does not change the formula from \fullref{action} except that we have to post-composing $R:\calO^{\bullet} \to R\calO^{\bullet}$.
On the other hand the $E^1$-term of Bousfield spectral sequence $E^r (R\calO^{\bullet} )$ for $\Tot R\calO^{\bullet}$
is equipped with the induced Gerstenhaber structure, whose formula is the same one from \fullref{turchin} with
$R$ post-composed. Thus again via $\varphi$ the Browder operation on $H_* (\Tot R\calO^{\bullet})$ and the Gerstenhaber
bracket on $E^{\infty}(R\calO^{\bullet})$ coincide with each other.\qed

%%% acknowledgement %%%

\subsection*{Acknowledgements}
The author expresses his great appreciation to Toshitake Kohno for his encouragement and advices.
The author is also grateful to Victor Turchin and James McClure for reading the draft of the previous version of
the paper and giving him many suggestions, to Dev Sinha for answering his questions, to Fred Cohen and
Ryan Budney for teaching him about the little disks actions and so on, and to Paolo Salvatore for kindly giving
the author his preprint.

This research is partially supported by the 21st century COE program at Graduate School of Mathematical Sciences,
the University of Tokyo.

\bibliographystyle{gtart}
\bibliography{link}

\end{document}